\documentclass{amsart}

\usepackage{amsmath, amsfonts, amssymb, amscd}

\newtheorem{theo}{Theorem}[section]
\newtheorem{lem}[theo]{Lemma}
\newtheorem{prop}[theo]{Proposition}

\newenvironment{pf}{\noindent{\it Proof.
}}{$\blacksquare$\par\medskip}

\newcommand{\C}{{\mathbb C}}
\renewcommand{\P}{{\mathbb P}}
\newcommand{\R}{{\mathbb R}}

\renewcommand{\H}{{\mathbb H}}

\newcommand{\g}{{\mathfrak g}}

\renewcommand{\u}{{\mathfrak u}}

\renewcommand{\a}{{\mathfrak a}}

\renewcommand{\t}{{\mathfrak t}}
\renewcommand{\b}{{\mathfrak b}}
\renewcommand{\c}{{\mathfrak c}}
\renewcommand{\d}{{\mathfrak d}}
\newcommand{\ee}{{\mathfrak e}}
\newcommand{\ff}{{\mathfrak f}}

\newcommand{\U}{{\operatorname {U}}}

\renewcommand{\varpi}{\psi}

\newcommand{\Gr}{{\operatorname{Gr}}}
\newcommand{\su}{{\mathfrak{su}}}
\newcommand{\so}{{\mathfrak {so}}}
\renewcommand{\sp}{{\mathfrak {sp}}}
\newcommand{\SU}{{\operatorname{SU}}}
\renewcommand{\SS}{{\operatorname{S}}}
\newcommand{\OO}{{\operatorname{O}}}
\newcommand{\SO}{{\operatorname{SO}}}
\newcommand{\Sp}{{\operatorname{Sp}}}
\newcommand{\Spin}{{\operatorname{Spin}}}
\newcommand{\G}{{\operatorname{G}}}
\newcommand{\SL}{{\operatorname{SL}}}
\newcommand{\CP}{{\C\!\operatorname{P}}}

\def\sideremark#1{\ifvmode\leavevmode\fi\vadjust{
\vbox to0pt{\hbox to 0pt{\hskip\hsize\hskip1em
\vbox{\hsize3cm\tiny\raggedright\pretolerance10000
\noindent #1\hfill}\hss}\vbox to8pt{\vfil}\vss}}}

\title[Homogeneous Lagrangian submanifolds of positive Euler characteristic ]
{Homogeneous Lagrangian submanifolds of positive Euler characteristic}
\author{Fabio Podest\`a}

\subjclass[2000]{53D12, 22F30} \keywords{Lagrangian submanifolds, homogeneous spaces}

\begin{document}

\begin{abstract} We fully classify all Lagrangian submanifolds of a complex Grassmannian which are an orbit of
a compact group of isometries and have positive Euler characteristic.\end{abstract}

\maketitle
\section{Introduction}
\bigskip

It is natural to try to find new examples of Lagrangian submanifold of K\"ahler manifolds among
orbits of compact Lie groups of holomorphic isometries. In \cite{BG}, the authors classified all the
Lagrangian submanifolds of the complex projective space which are homogeneous under the action of a compact
simple Lie group $G$ of isometries, while apparently there exist many more examples when one drops the
assumption of simplicity of $G$; interestingly, not all the examples found in this
classification have parallel second fundamental form and therefore they provide new examples. More recently,
H. Ma and Y. Ohnita (\cite{Om}) classified homogeneous Lagrangian submanifolds of the
complex quadrics via a neat correspondence with homogeneous isoparametric hypersurfaces of the standard sphere.\par
However, a full classification of homogeneous Lagrangian submanifolds (HLS) of Hermitian symmetric spaces or even of
the complex projective space seems
out of reach by now. In this paper we focus on a particular class of HLS of the complex Grassmannian $\Gr_k(\C^n)$,
($1\leq k\leq [\frac{n}{2}]$) endowed with the standard K\"ahler structure,
namely those which have positive Euler characteristic; in this case we are able to provide a full classification.
\begin{theo}\label{main} Let $G$ be a compact connected Lie subgroup of $\U(n)$ acting almost faithfully on the
complex Grassmannian $\Gr_k(\C^n)$ ($1\leq k\leq [\frac{n}{2}]$) with a Lagrangian
orbit $\mathcal L$  with positive Euler characteristic.
Then \par
\noindent (a)\ $G$ is simple, except for $G=\SO(4)$ and $\mathcal L = \Gr_2(\R^4)\subset \Gr_2(\C^4)$;\par
\noindent (b)\  the
orbit $\mathcal L$ is the only isotropic $G$-orbit and it is totally geodesic; moreover
the triple $(G,\mathcal L,\Gr_k(\C^n))$ appears in Table 1,
where $\Gr_p(\R^n)$ and $\Gr_p(\H^n)$ denote the Grassmannian of real $p$-planes
in $\mathbb R^n$ and the Grassmannian of quaternionic $p$-planes in $\H^n$ respectively. \par
In particular $\R P^{2n}$ is the only HLS with positive
Euler characteristic in a complex projective space.\par
Conversely, every subgroup $G\subset \U(n)$ appearing in Table 1 has $\mathcal L$ as a Lagrangian orbit of
positive Euler characteristic.\medskip
\begin{table}[h]
\begin{tabular}{|c|c|c|}


\hline $G$ & $\mathcal L$ & $\Gr_k(\C^n)$ \\
\hline $\SO(n)$ &$\SO(n)/\SS(\OO(p)\times\OO(n-p))\cong \Gr_{p}(\R^n)$,\ $p(n-p)$ $\operatorname{even}$ & $\Gr_{p}(\C^n)$
\\
\hline $\Sp(m)$ &$\Sp(m)/\Sp(p)\times\Sp(m-p)\cong \Gr_{p}(\H^m)$ & $\Gr_{2p}(\C^{2m})$
\\
\hline $\Spin(7)$ &$\SO(7)/\U(3)\cdot\mathbb Z_2\cong \Gr_2(\R^8)$ & $\Gr_2(\C^8)$
\\
\hline $\operatorname{G}_2$ &$\operatorname{G}_2/\U(2)\cdot\mathbb Z_2\cong\Gr_2(\R^7)$ & $\Gr_2(\C^7)$\\
\hline
\end{tabular}\vspace{0.1em}\medskip
\caption{}
\end{table}\par
\end{theo}
We remark here that the uniqueness of a homogeneous Lagrangian $G$-orbit, when
$G$ is semisimple, is a general fact (see \cite{BG}). \par
In the next section, we will set up notations and give the proof of Theorem 1.1; the main tools used in the
proof are the moment maps and the representation theory for compact Lie groups. We also remark that
the complexified group $G^\C$
acts on the Grassmannian with an open orbit if $G$ has a Lagrangian orbit (see the proof of Lemma \ref{uni}) , so that
another approach could be to use the classification of prehomogeneous irreducible
vector spaces by Sato and Kimura (\cite{SK}); we prefer to give an independent proof together with a more theoretical
view of the geometric setting. \par
There are other situations, where a full classification of Lagrangian orbits is feasible; for instance, the following observation
can be easily deduced from the classification by Brion of multiplicity free homogeneous spaces (\cite{Br})
\begin{prop} Let $G$ be a compact Lie group acting isometrically and holomorphically on a K\"ahler manifold $M$ in
a Hamiltonian fashion. If the action of
$G$ is multiplicity free and $G$ has a Lagrangian orbit $\mathcal L$ of positive Euler characteristic, then $G$ is semisimple
and $\mathcal L$ is locally diffeomorphic to the product of symmetric spaces of inner type.\end{prop}
A complete list of such actions can be easily obtained using \cite{Br}.

\bigskip
\bigskip
\section{Preliminaries and proof of the main theorem}
\bigskip
Throughout the following $G$ will denote a compact connected Lie group and $T\subset G$ a fixed maximal torus; gothic letters
will indicate the corresponding Lie algebras. We also recall that a $G$-homogeneous space $G/H$ has positive
Euler characteristic if and only if the subgroup $H$ has maximal rank in $G$ (see e.g. \cite{On}).\par
If $G$ is a subgroup of $\U(n)$ for some $n\in \mathbb N$, then
$G$ acts isometrically on every complex Grassmannian $\Gr_k(\C^n)$ ($1\leq k\leq [\frac{n}{2}]$) endowed with its standard
K\"ahlerian structure; we will suppose that the action is almost faithful, i.e. the subgroup of $G$ given by all
elements of $G$ that act trivially on the
Grassmannian is a finite group. Moreover, the $G$-action
is Hamiltonian, i.e. there exists a moment map $\mu:\Gr_k(\C^n)\to \g^*$, uniquely determined up
to translation, which is $G$-equivariant and satisfies for $p\in \Gr_k(\C^n)$, $Y\in T_p\Gr_k(\C^n)$ and $X\in \g$
$$\langle d\mu_p(Y),X\rangle = \omega_p(Y,\hat X),$$
where $\omega$ is the K\"ahler form and $\hat X$ denotes the Killing vector field induced by $X$. We recall that
in general an orbit $\mathcal O$ is isotropic if and only if $\mu(\mathcal O)$ is contained in the annihilator
$[\g,\g]^o := \{\phi\in \g^*|\ \langle\phi,[\g,\g]\rangle = 0\}$, in particular if and only if $\mu(\mathcal O) = 0$
whenever $\g$ is semisimple (see e.g. \cite{Ki}). \par
Now, if we embed $\Gr_k(\C^n)$ into ${\mathbb P}(\Lambda^k(\C^n))$ via the Pl\"ucker embedding, we may use the
expression for the moment mapping of a Lie subgroup of the group of all isometries of a complex projective space (see e.g.
\cite{Ki}):
in this way if $X\in \g$ and $v_1,\ldots,v_k$ is a basis of a $k$-plane $\pi\in \Gr_k(\C^n)$, then there is a
moment map $\mu:
\Gr_k(\C^n)\to \g^*$ such that
$$\langle \mu(\pi),X\rangle = \frac{1}{||v_1\wedge\ldots\wedge v_k||^2}(X\cdot v_1\wedge v_2\wedge\ldots\wedge v_k +
\ldots$$
$$\ldots +\ v_1\wedge\ldots\wedge X\cdot v_k,\ v_1\wedge\ldots\wedge v_k), \eqno (1.1)$$
where $(\cdot,\cdot)$ denotes the standard Hermitian product in $\Lambda^k(\C^n)$ induced by the standard $\U(n)$-invariant
Hermitian product on $\C^n$. We now suppose to have a $G$-orbit $\mathcal L = G\cdot \pi = G/H$, with
$\pi\in \Gr_k(\C^n)$, of positive Euler characteristic; then we can suppose
$T\subseteq H$. We can then consider the weight space decomposition $\C^n = \bigoplus_{\lambda\in \Lambda}V_\lambda$
w.r.t. to the torus $T$, where $\Lambda\subset \operatorname{Hom}(\t, i\mathbb R)$ is the set of weights; the $k$-plane
$\pi$, being $T$-invariant, splits as
$\pi = \bigoplus_{\lambda\in \Lambda}\pi_\lambda$, where $\pi_\lambda := \pi\cap V_\lambda$.\par
 If $\mathcal L$ is Lagrangian,
the isotropy representation of $H$ at $\pi$ is isomorphic to $\imath^\C$, where $\imath$ is the isotropy
representation $\imath: H\to O(T_\pi\mathcal L)$. We now claim that
$$\operatorname{for\ every } \lambda\in \Lambda,\quad\qquad V_\lambda\cap\pi
\not=\{0\}\Rightarrow V_\lambda\subseteq \pi.$$
Indeed suppose that there exists $\lambda$ with $V_\lambda\cap \pi\neq\{0\}$ but $V_\lambda$ not fully contained in $\pi$; then
$V_\lambda$ intersects the orthogonal $\pi^\perp$ non trivially. Since the tangent space $T_\pi\Gr_k(\C^n) \cong
\operatorname{Hom}(\pi,\pi^\perp)\cong \pi^*\otimes\pi^\perp$, we see that
the isotropy representation of $H$ at $\pi$ has $-\lambda+\lambda = 0$ as a weight; but the isotropy
representation of $H$, given by $\imath^\C$, does not have $0$ as a weight because $H$ contains a maximal torus, a contradiction.
So there is a subset $Q\subset \Lambda$ such that
$$\pi = \bigoplus_{\lambda \in Q}V_\lambda.$$
Moreover, since the weights of $\imath^\C$ are roots of $\g^\C$ w.r.t. the Cartan subalgebra $\t^\C$, we see that
the weights of $\pi^*\otimes\pi^\perp$ are roots, i.e.
$$ \mu-\lambda \in R,\quad \forall \lambda\in Q,\ \forall \mu\in \Lambda\setminus Q\eqno (1.2)$$
where $R$ is the root system of $\g^\C$ relative to $\t^\C$.\par
\begin{lem}\label{uni} If $G$ has a Lagrangian orbit $\mathcal L$ in $\Gr_k(\C^n)$ with positive Euler characteristic,
then $G$ is semisimple and it acts on $\C^n$ irreducibly. Moreover the orbit $\mathcal L$ is the only isotropic orbit in
$\Gr_k(\C^n)$.
\end{lem}
\begin{pf} \noindent (1)\ We prove that $G$ is semisimple. Indeed we claim that the center $Z$ of $G$ is a finite
group; in order to prove this, we will show that $Z$ acts trivially on the Grassmannian and our claim
will follow form the fact that the $G$-action is almost faithful. Now, $Z$ acts trivially on the orbit
$\mathcal L = G\cdot p = G/H$ because $H$ has maximal rank and therefore contains $Z$; so the isotropy
representation of every $z\in Z$ is trivial on the tangent space $T_p\mathcal L$, hence on $JT_p\mathcal L$ , where $J$ is the complex
structure. Since $T_p\mathcal L$ and $JT_p\mathcal L$ span the whole tangent space at $p$, $z$ acts as the identity. \par
\noindent (2)\ We now prove that $G$ acts irreducibly. If not, we split
$\C^n = W_1\oplus W_2$ for some non-trivial $G$-invariant complex subspaces $W_1,W_2$; since
$\Lambda = \Lambda_1\cup \Lambda_2$
where $\Lambda_i$ is the weight system of $W_i$ ($i=1,2$), we have that $\pi = (\pi\cap W_1)\oplus (\pi\cap W_2)$. On
the other hand if $\mathcal L$ is Lagrangian, the complexified group $G^\C$ acts on the Grassmannian with the open
orbit $G^\C\cdot \pi$: indeed, if $J$ denotes the complex structure,
$T_\pi(G^\C\cdot\pi) = \{\hat X_\pi + J\hat Y_\pi|\ X,Y\in \g\} = T_\pi\mathcal L + JT_\pi\mathcal L = T_\pi
\Gr_k(\C^n)$. But any $k$-plane $\pi'= u\cdot \pi$ ($u\in G^\C$) in the $G^\C$-orbit will be of the form
$\pi' = (\pi'\cap W_1)\oplus (\pi'\cap W_2)$ because $G^\C$ preserves the decomposition $\C^n = W_1\oplus W_2$;
since any $k$-plane of the form $\pi'$ is not generic, the orbit $G^\C\cdot \pi$ is not open, a contradiction.\par
\noindent (3)\ Since $G$ is semisimple, an orbit $\mathcal O$ is isotropic if and only if $\mathcal O \subseteq \mu^{-1}(0)$;
we now claim that $\mu^{-1}(0)$ coincides with the Lagrangian orbit $\mathcal L$, proving our last assertion. We first
recall that the set $\mu^{-1}(0)$ is connected (see \cite{Ki}), so that it is will be enough to show that $\mathcal L$
is open in $\mu^{-1}(0)$. We recall that the set $G^\C\cdot \mathcal L$ is open in $\Gr_k(\C^n)$ because $\mathcal L$ is
Lagrangian; we claim that $(G^\C\cdot \mathcal L)\cap \mu^{-1}(0) = \mathcal L$, showing that $\mathcal L$ is open
in $\mu^{-1}(0)$. This last assertion follows from Kirwan's Lemma (\cite {Ki}, p. 97): if two points in $\mu^{-1}(0)$ lie in the same
$G^\C$-orbit, they belong to the same $G$-orbit.  \end{pf}
Since $G$ is semisimple, the orbit $G\cdot\pi$ is isotropic if and only if $\mu(\pi)=0$; we claim that this condition is
equivalent to
$$\sum_{\lambda \in Q} m_\lambda \cdot \lambda = 0,\eqno (1.3)$$
where $m_\lambda$ is the multiplicity of the weight $\lambda$, i.e. $m_\lambda = \dim V_\lambda$; indeed, if we
choose a basis of $\pi$ made of
weight vectors, then (1.1) shows that (1.3) holds if and only if $\langle \mu(\pi),H\rangle = 0$ for every
$H\in \t$. A direct inspection shows that for any $E$ in any root space $\g_\alpha\subset \g^\C$
we have that $(E\cdot\pi,\pi)=0$ and we get our claim.\par
Moreover, since the $H$-representation $\pi^*\otimes\pi^\perp$ is equivalent to $\imath^\C$, that has root
spaces as weight spaces, we see that $\dim V_\mu\otimes V_\lambda^* = 1$ for every
$\lambda\in Q$ and $\mu\in \Lambda\setminus Q$, i.e.
$$m_\lambda = 1, \qquad \forall\ \lambda \in \Lambda.\eqno (1.4)$$

\noindent{\bf Remark.}\ Using these considerations and notations, we can easily provide examples
of {\it isotropic\/} orbits of positive Euler characteristic
in a Grassmannian $\Gr_k(\C^n)$. Indeed the $G$-orbit through the plane $\pi = \bigoplus _{w\in W} V_{w\cdot \lambda_o}$,
where $W$ is the Weyl group of $G$ and $\lambda_o$ is any weight whose $W$-orbit does not exhaust the set of all
weights $\Lambda$ is isotropic (because $\sum_{w\in W}w(\lambda_o) = 0$) and the stabilizer $G_\pi$ has maximal rank in
$G$; however these orbits are rarely Lagrangian.\par
\medskip
We now consider the case $k=1$, i.e. the complex projective space $\CP^{n-1}$; our result can be easily deduced from
Lemma 2.1 and the list in \cite{BG} or the tables in \cite{SK}, but we prefer to give here a more direct proof.\par
\begin{lem}\label{simple1} If $k=1$, then $G$ is simple.\end{lem}
\begin{pf} Suppose $G$ is (locally) isomorphic to a product of simple groups
$G_1\times\ldots\times G_s$
for some $s\geq 2$ and accordingly $V = \bigotimes_{i=1}^s W_i$ for some irreducible $G_i$-modules $W_i$, $i=1,\ldots,s$;
the maximal torus $T$ splits (locally) as $T_1\times\ldots\times T_s$ for some maximal tori $T_i\subset G_i$ and the
weight system $\Lambda = \Lambda_1+\dots+\Lambda_s$, where $\Lambda_i$ is the weight system (relative to $T_i$) of
$G_i$ acting on $W_i$, $i=1,\ldots,s$. The $1$-plane $\pi$ is given by some weight space
$V_\lambda = \bigotimes W_{i,\lambda_i}$ with $\lambda = \lambda_1+\ldots+\lambda_s$
for some weights $\lambda_i\in \Lambda_i$, ($i=1,\ldots,s$); if we now select $\mu_i\in \Lambda_i$ with $\mu_i\neq \lambda_i$,
we see that the weight $\mu := \sum_{i=1}^s\mu_i\neq \lambda$ and $\mu-\lambda$ is not a root , contradicting (1.2).
So $G$ is simple.\end{pf}
\par
We now conclude the case $k=1$. By (1.3) we have that $0\in \Lambda$ and $\pi = V_0$; moreover by (1.2) we have $\Lambda\subseteq R$.
This implies that the representation is the adjoint representation of a Lie algebra of
rank $1$ or $\g$ is of type $\b_{l}, \c_{l}, \ff_4, \g_2$ and the irreducible
representation has highest weight $\lambda_1,\lambda_2,\lambda_4,\lambda_1$ respectively,
where $\lambda_i$ denote the fundamental weights (see e.g. \cite{He}, p. 523). Among these, only the case $\g=\so(2l+1)$
acting on $\C^{2l+1}$ in the standard way or $\g=\g_2$ acting on $\C^7$ in the standard way are admissible: indeed
when $\g=\c_l$ ($l\geq 3$) is acting on $\Lambda_o^2(\C^{2l})$ one easily sees that the zero weight space
has dimension $l-1$ (or also it has no open orbits in the corresponding projective space by \cite{SK}, p.105); when
$\g=\ff_4$ is acting on $\C^{26}$ in the standard way, the zero weight space is two dimensional (and also the action
on $\CP^{25}$ has no open orbits by \cite{SK} p. 143) and therefore it can be ruled out. \par

We now deal with the case $k\geq 2$.
\begin{lem} If $k\geq 2$ then $G$ is simple, unless $G=\SO(4)$ acting on $\Gr_2(\C^4)$ with $\mathcal L =
\Gr_2(\R^4)$ as a Lagrangian orbit.\end{lem}
\begin{pf} Suppose $G$ is not simple and use the same notations as in Lemma \ref{simple1}, with $n_i := \dim W_i$
and $n_1\geq n_2\geq\ldots\geq n_s\geq 2$. First suppose $s\geq 3$:  If $\lambda = \lambda_1+\ldots
+\lambda_s$ is in $Q$, then by (1.2), $Q$ contains the disjoint subsets
$(\Lambda_1\setminus\{\lambda_1\})+(\Lambda_2\setminus\{\lambda_2\})+
\Lambda_3+\ldots+\Lambda_s$, $\lambda_1+(\Lambda_2\setminus\{\lambda_2\})+
(\Lambda_3\setminus\{\lambda_3\})+\Lambda_4+\ldots+\Lambda_s$ and
$(\Lambda_1\setminus\{\lambda_1\}) + \lambda_2 +
(\Lambda_3\setminus\{\lambda_3\})+\Lambda_4+\ldots+\Lambda_s$; therefore
$k\geq [(n_1-1)(n_2-1)n_3+(n_2-1)(n_3-1)+(n_1-1)(n_3-1)]n_4\ldots n_s$. Since
we suppose $k\leq \dim V/2$, we have that $n_1n_2n_3\leq 2(n_1+n_2+n_3)-4$; from this we easily
get that $n_i=2$ for all $i=1,\ldots,s$. Therefore $G_i\cong \SU(2)$ and $\dim \mathcal L = 2s = k(2^s-k) \geq 2(2^s-2)$,
i.e. $s\geq 2^s-2$ and this contradicts $s\geq 3$.\par
When $s=2$, the subset $Q$ contains $(\Lambda_1\setminus\{\lambda_1\})+(\Lambda_2\setminus\{\lambda_2\})$, hence
$n_1n_2/2\geq k\geq (n_1-1)(n_2-1) + 1$. This implies $n_2=2$ and $k\geq n_1$, hence $G\cong G_1\times \SU(2)$,
$V=W_1\otimes \C^2$ and $\dim\mathcal L = \dim G_1/H_1 + 2$ for some maximal rank subgroup $H_1\subset G_1$. Therefore
$\dim G_1/H_1 + 2 = k(2n_1-k) \geq n_1^2$; if the rank of $G_1$ is greater than one, then $\dim H_1\geq 2$ and
$\dim G_1\geq n_1^2$: this never happens because $G_1$ is simple, hence $G_1$ is locally isomorphic to a subgroup
of $\SL(n_1)$ and therefore $\dim G_1\leq n_1^2-1$. \par
We are left with the case $\operatorname{rank}(G_1)=1$, namely $G_1\cong \SU(2)$ and $W\cong \C^2\otimes\C^2$;
this gives the Lagrangian orbit $\Gr_2(\R^4)$ inside the Grassmannian $\Gr_2(\C^4)$. \end{pf}

The following Lemma will be very useful (see also \cite{BG})
\begin{lem}\label{lem} Let $G$ be a (semi)simple compact Lie group acting linearly on some complex vector space $V$ and
therefore on
the Grassmannian $\Gr_k(V)$ for $1\leq k\leq [\frac{n}{2}]$, with $\dim V = n$. Suppose that the representation of $G$
on $V$ is of real (quaternionic and $k$ is even) type.
Then $G$ has a Lagrangian orbit in $\Gr_k(V)$ if and only if
it acts transitively on a real (quaternionic resp.) Grassmannian,
i.e on $\SO(n)/\SO(k)\times \SO(n-k)$ ( $\Sp(n/2)/\Sp(k/2)\times\Sp((n-k)/2)$ resp.).
\end{lem}
\begin{pf} We denote by $\mu:\Gr_k(V)\to \g^*$ the unique moment map for the $G$-action on the Grassmannian; it is
well known that every isotropic $G$-orbit lies in the set $\mu^{-1}(0)$, which is connected (see \cite{Ki}). If
$G$ has a Lagrangian orbit $\mathcal L$, then it admits no other isotropic orbit by Lemma \ref{uni}. If the representation is of real type (quaternionic type),
then $G\subset \SO(n)$ ( $G\subset \Sp(n/2)$ resp.) and $\SO(n)$ ($\Sp(n/2)$ resp.) does have a Lagrangian orbit $\mathcal L$; if $G$ has
a Lagrangian orbit $\mathcal L'$, it also has an isotropic orbit inside $\mathcal L$, hence $\mathcal L = \mathcal L'$.
\end{pf}
\par
Since $G^\C$ has an open orbit in the Grassmannian and the generic stabilizer has full rank, we have that
$$\dim G - \operatorname{rank}(G) \geq k(n-k) \geq 2(n-2),\eqno (1.5)$$
hence
$$n\leq \frac{1}{2} (\dim G - \operatorname{rank}(G) + 4).\eqno (1.6)$$
Using the well known classification of irreducible
representation of simple groups whose representation space $V$
satisfies $\dim V\leq \dim G$ (see e.g. \cite {FH}, p.414), we easily see that the condition (1.6)
is satisfied only for the groups and representations listed in Table 2 (where we use the
standard notation for fundamental weights and where we do not list dual representations of admissible ones)
\medskip
\begin{table}[!h]
\begin{tabular}{|c|c|c|c|}

\hline $\g$ & highest weight & $\frac{1}{2} (\dim \g - \operatorname{rank}(G) + 4)$ & $\dim V$\\
\hline $\a_n$ &$\lambda_1,\lambda_2$ & $(n^2+n+4)/2$&
$n+1, \frac{1}{2} n(n+1)$\\
\hline $\b_n,\ n\geq 2$ &$\lambda_1, \lambda_n (n=2,3,4)$ & $n^2+2$&
$2n+1,4,8,16$\\
\hline $\c_n,\ n\geq 3$ &$\lambda_1$ & $n^2+2$&
$2n$\\
\hline $\d_n,\ n\geq 4$ &$\lambda_1, \lambda_{n-1} (n=4,5,6)$ & $n^2-n+2$&
$2n,8,16,32$\\
\hline $\ee_6$ &$\lambda_1$ & $38$&
$27$\\
\hline $\ee_7$ &$\lambda_1$ & $65$&
$56$\\
\hline $\ee_8$ &$-$ & $122$&
$-$\\
\hline $\ff_4$ &$\lambda_4$ & $26$&
$26$\\
\hline $\g_2$ &$\lambda_1$ & $8$&
$7$\\
\hline
\end{tabular}\vspace{1em}
\caption{}\label{Tten}
\end{table}\par
We now study the case $k=2$. When $\g$ is of classical type, we will first discuss only the representations in the table
different from $\lambda_1$ and we will postpone the discussion of these to the end of the proof.\par
If $\g$ is of type $\a_n$, then $\lambda_2$ is not possible unless $n=3$: indeed the $2$-plane
$\pi$ must be spanned by two opposite weights and this representation has never opposite weights unless $n=3$ and
$\g = \su(4)\cong\so(6)$ is acting on $\C^6$.\par
When $\g$ is of type $\b_n$, we have to discuss different special cases, namely: (1) $n=2$ and the representation is the
spin representation : this is admissible, because $\so(5)\cong\sp(2)$ and the spin representation goes over to the standard one;
(2) $n=3$ and $\so(7)$ acts on $\C^8$ via spin : the group $\Spin(7)$ acts transitively on the real Grassmannian
$\SO(8)/\SO(2)\times\SO(6)\cong\SO(8)/\U(4) = \SO(7)/\U(3)$, hence on the Lagrangian orbit $\SO(8)/\SS(\OO(2)\times\OO(6))$
which can be written as $\SO(7)/\U(3)\cdot\mathbb Z_2$;\ (3)\ $n=4$ and $\so(9)$ acts on $\C^{16}$ via spin: it
is a real representation and $\Spin(9)\subset\SO(16)$; since $\Spin(9)$ does not act transitively on the real two plane
Grassmannian (see \cite{On}), we rule this out by Lemma~\ref{lem}.
\par
If $\g$ is of type $\d_n$ we just have to discuss the half-spin representations listed in the table: (1) for $n=4$
we have $\so(8)$ acting on $\C^8$ and this is admissible; (2) for $n=5$ we have $\so(10)$ acting on $\C^{16}$ via
half-spin. In either
case, the weights lie in a single Weyl orbit and we can suppose that the $2$-plane $\pi$ contains the highest
weight $\lambda$; but $-\lambda$ is not a weight and we can rule this case out. (3) $n=6$ and $\so(12)$ is acting on
$\C^{32}$ via half-spin. Again, we may suppose that $\pi$ is generated by the weight spaces of $\lambda, -\lambda$,
where $\lambda$ is the highest weight; so the action has an isotropic orbit in $\Gr_2(\C^{32})$. Since
$\g_\pi$ contains $\g_\lambda := \{X\in \g;\ X\cdot V_\lambda\subseteq V_\lambda\}\cong \u(6)$ and since $\u(6)$ is maximal
in $\so(12)$, we see that $\g_\pi\cong\u(6)$ and $\dim_\R G\cdot\pi= 30$, while $\dim_\C\Gr_2(\C^{32}) = 60$; so this
case is impossible.\par
When $\g=\ee_6$, the two representations $\lambda_1$ and $\lambda_6$ are dual to each other and their weight systems are
$\Lambda_1 = -\Lambda_6$; moreover each of them consists of a single Weyl orbit. Therefore the $2$-plane $\pi$ can be
supposed to contain the highest weight space (using the action of the Weyl group); on the other it should contain
the opposite of the highest weight, that is not a weight. \par
When $\g=\ee_7$, the irreducible representation $\lambda_1$ is of quaternionic type and has complex dimension $56$; the group $\operatorname{E}_7$
does not act transitively on $\H\P^{27}$ (see \cite{On}), hence this case can be ruled out by Lemma~\ref{lem}.
\par
For $\g=\ff_4$, it is well known that $26$-dimensional irreducible representation is of real type, while the group
$\operatorname F_4$ does not act transitively on any Grassmannian (see \cite{On}); so again we rule this out
by Lemma~\ref{lem}. \par
For $\g=\g_2$ acting on $\C^7$, we recall that the group $G_2$ is contained in $\SO(7)$ and that it acts transitively
on the real Grassmannian $\SO(7)/\SS(\OO(2)\times\OO(5))\cong \G_2/\U(2)\cdot\mathbb Z_2$ (see e.g. \cite{On}) which is
Lagrangian in $\Gr_2(\C^7)$.\par
We now focus on the higher dimensional Grassmannians, i.e. $k\geq 3$. The inequality (1.5) becomes
$$n \leq \frac{1}{3}(\dim G - \operatorname{rank}(G) + 9) \eqno (1.7)$$
and a direct inspection shows that the only possible representation in Table 2 are listed in Table 3.
\medskip
\begin{table}[!h]
\begin{tabular}{|c|c|c|c|}

\hline $\g$ & highest weight & $\frac{1}{3} (\dim \g - \operatorname{rank}(G) + 9)$ & $\dim V$\\
\hline $\a_n$ &$\lambda_1,\lambda_2 (n\leq 3)$ & $(n^2+n+9)/3$&
$n+1, \frac{1}{2} n(n+1)$\\
\hline $\b_n,\ n\geq 2$ &$\lambda_1, \lambda_n (n=2,3)$ & $(2n^2+9)/3$&
$2n+1,4,8$\\
\hline $\c_n,\ n\geq 3$ &$\lambda_1$ & $(2n^2+9)/3$&
$2n$\\
\hline $\d_n,\ n\geq 4$ &$\lambda_1, \lambda_{n-1} (n=4,5)$ & $(2n^2-2n+9)/3$&
$2n,8,16$\\
\hline $\ee_6$ &$\lambda_1$ & $27$&
$27$\\
\hline $\ee_7$ &$-$ & $45$&
$-$\\
\hline $\ee_8$ &$-$ & $83$&
$-$\\
\hline $\ff_4$ &$-$ & $19$&
$-$\\
\hline $\g_2$ &$\lambda_1$ & $7$&
$7$\\
\hline
\end{tabular}\vspace{1em}
\caption{}\label{Tten}
\end{table}\par
Now, for $\g=\a_n$, we see that the second fundamental representation for $n\leq 3$ is admissible and corresponds to
the dual of $\lambda_1$ for $\SU(3)$ and to $\lambda_1$ of $\so(6)$ for $n=3$. For $\g=\b_n$, $\lambda_2$ with $n=2$
is the first fundamental of $\sp(2)$, while we have to discuss $\lambda_3$ for $n=3$; since $\Spin(7)\subset \SO(8)$ and
$\SO(8)$ has a Lagrangian orbit in $\Gr_3(\C^8)$ given by a real Grassmannian on which $\Spin(7)$ does not act
transitively (see e.g. \cite{On}), this case is ruled out by Lemma~\ref{lem}.\par
If $\g=\d_n$, we just have to discuss the case $\lambda_5$ for $n=5$, i.e. $\Spin(10)$ acting on $\C^{16}$ via 
half-spin. A Lagrangian orbit in $\Gr_3(\C^{16})$ would have dimension $39$, while there is no subgroup of $\Spin(10)$ of 
maximal rank ($=5$) and dimension $\dim\Spin(10)-39 = 6$: indeed such a subgroup $H$ would contain a maximal torus $T$ 
and $\dim H/T=1$, forcing $H$ to be abelian, a contradiction.\par
If $\g=\ee_6$ the representation space is $\C^{27}$ and all the weights lie in a single Weyl orbit (see \cite{lie}); if it
has a Lagrangian orbit in $\Gr_3(\C^{27})$, we can suppose that the $3$-plane $\pi\in \Gr_3(\C^{27})$ contains the
highest weight space. Now, $\pi^\perp$ is the
sum of $24$ one-dimensional weight spaces $V_\mu$, where $\mu = \lambda_1 -\alpha_\mu$, $\alpha_\mu$ being a positive root
with $\langle\lambda_1,\alpha_\mu\rangle\not= 0$. So there would be at least $24$ positive roots whose expression in terms of
simple roots $\{\alpha_1,\ldots,\alpha_6\}$ have a non zero coefficient relative to $\alpha_1$; but there are only $20$ such roots
in $\ee_6$ (see e.g. \cite{He}). \par
The case $\g_2$ can be also ruled out using Lemma~\ref{lem} and the fact that $\operatorname{G}_2$ does not act transitively
on $\Gr_3(\mathbb R^7)$.\par
Using Table 2 it is immediate to see that when $k\geq 4$, there are no other admissible representation but
the first fundamental ones.\par
We now have to discuss the first fundamental representations; for $\b_n$ and $\d_n$, we know that the real Grassmannians are
Lagrangian orbits and therefore there are no others; we just have to discuss the case of $\g=\c_n$ acting on a
Grassmannian $\Gr_{2k+1}(\C^{2n})$ for some $1\leq k\leq n-1$. In this case, there is no Lagrangian orbit: indeed,
suppose the orbit through $\pi\in\Gr_{2k+1}(\C^{2n})$ is Lagrangian; now, $\pi$ would be spanned by an odd number of
weight spaces relative to some weights $\mu_1\ldots\mu_{2k+1}$ with $\sum_{i=1}^{2k+1}\mu_i = 0$. But
the root system $R$ of $\c_n$ is given by $R = \{\pm 2e_i\ (1\leq i\leq n),\ \pm e_i\pm e_i\ (1\leq i\not= j\leq n)\}$
(see e.g. \cite{He} p. 463) and the weight system is simply $\Lambda = \{\pm e_i,\ 1\leq i\leq n\}$; so it is immediate
to see that the condition $\sum_{i=1}^{2k+1}\mu_i = 0$ cannot happen.

\bigskip
\bigskip
\font\smallsmc = cmcsc8 \font\smalltt = cmtt8 \font\smallit = cmti8
\hbox{\parindent=0pt\parskip=0pt \vbox{\baselineskip 9.5 pt
\hsize=3.1truein \obeylines {\smallsmc Fabio Podest\`a
Dip. Matematica e Appl. per l'Architettura
Universit\`a di Firenze
Piazza Ghiberti 27
I-50100 Firenze
ITALY }\medskip
{\smallit E-mail}\/: {\smalltt podesta@math.unifi.it } } }

\end{document}